\documentclass[12pt]{amsart}
\usepackage{amsthm,amsmath,amssymb}
\usepackage{addfont}
\usepackage{mathrsfs}
\usepackage[T1]{fontenc}
\numberwithin{equation}{section}

\def\ca{{\mathcal A}}

\def\cas{{\mathcal S}}

\def\bc{{\mathbb C}}

\def\bn{{\mathbb N}}

\def\g{\gamma}

\def\l{\lambda}

\def\f{\varphi}

\theoremstyle{plain}
\newtheorem{lemma}{Lemma}[section]
\newtheorem{proposition}[lemma]{Proposition}
\newtheorem{theorem}[lemma]{Theorem}
\newtheorem{corollary}[lemma]{Corollary}
\theoremstyle{definition}

\newtheorem{remark}[lemma]{Remark}
\newtheorem{definition}[lemma]{Definition}
\newtheorem{example}[lemma]{Example}

\begin{document}

\title[ Uniqueness of some matrix factorizations ]{\textsc{  Unique Matrix Factorizations associated to bilinear forms and Schur multipliers  }}

\author[E.~Christensen]{Erik Christensen}
\address{\hskip-\parindent
Erik Christensen, Mathematics Institute, University of Copenhagen, Copenhagen, Denmark.}
\email{echris@math.ku.dk}
\date{\today}

\begin{abstract}
Grothendieck's inequalities for operators and bilinear forms imply some factorization results for complex  $m \times n$  matrices. The theory of operator spaces provides a set up which describes 4 norm optimal factorizations of Grothendieck's sort. It is shown  that 3 of the optimal factorizations are uniquely determined  and the remaining one is  unique in some cases.
 \end{abstract}

\maketitle

\section{Introduction and Notation}

The content of this article is closely related to our recent article \cite{C3}, 
 {\em Bilinear forms, Schur multipliers, complete boundedness and duality}, where we studied a complex   $m \times n $  matrix $X$ from 4 different points of views, and then applied the operator space theory and Grothendieck's insights to the matrices we look at. We quote from \cite{C3}. 
  
{\em $"$A complex  $m \times n $ matrix $X$ may represent many different things in pure and applied mathematics. In this article we will focus on the interpretations of $X$ in  4 different ways \begin{itemize} 
\item[(i)] 
As the matrix for a linear mapping $F_X$ of the $n$  dimensional abelian  C*-algebra $\ca_n \,:= \,C(\{1, \dots, n\}, \bc)$ into the $m$ dimensional Hilbert space $\bc^m.$
\item[(ii)]
As the kernel for a bilinear form $B_X$ on the product $\ca_m \times \ca_n$ of C*-algebras given by $$B_X(a, b) \, := \, \sum_{i = 1}^m\sum_{j=1}^n X_{(i,j)}a(i)b(j) .$$
\item[(iii)] 
 As a  a linear mapping $S_X$ on $(M_{( m \times n)}(\bc), \|.\|_\infty)$ induced by Schur multiplication by $X$  - or entry wise multiplication - given by $$S_X(A)_{(i,j)}\, : = \, X_{(i,j)}A_{(i,j)}.$$ 
\item[(iv)] 
  As a  a bilinear mapping $T_X$ of $\ca_m \times M_{( m \times n)}(\bc)$ into the Hilbert space $\bc^n$ given as   $$T_X(a,B)_j \, : = \, \sum_{i =1}^m a(i) X_{(i,j)}B_{(i,j)}."$$
\end{itemize}}

 The main results of our previous article \cite{C3} are expressed in terms of the completely bounded norms of the 4 objects  $F_X, B_X, S_X$ and $T_X.$ We will not repeat the descriptions from \cite{C3} regarding {\em complete boundedness} and {\em operator space theory} but mention that the article \cite{C3} contains a short introduction to the most needed results from that theory, and the  text books \cite{ER} by Effros and Ruan, \cite{Pa} by Paulsen and \cite{Pi2} by Pisier contain much more than needed here. The suffix $cb$ will appaer in several places and it stands for the words {\em completely bounded. } This concept is easy to define so we will repeat the definition here.
 \begin{definition} Let $H$ and $K$ be complex Hilbert spaces, $\cas$ a subspace of $B(H)$ and $\f : \cas \to B(K)$ a bounded  linear map. For each natural number $n$ we define $\f_n : \cas \otimes M_n(\bc) \to B(K) \otimes M_n(\bc) $ by $\f_n := \f \otimes \mathrm{id}_{M_n(\bc)} .$ If the set of norms $\{\|\f_n\|\, : \, n \in \bn\,\}$ is bounded, we say that $\f$ is completely bounded and define its completely bounded norm $\|\f\|_{cb}$ as the supremum of this set of norms.
  \end{definition} 
 On the other hand the results in \cite{C3} show that certain norms of completely bounded linear  or bilinear maps may be expressed as optimal values for some factorization properties of a matrix. With this aspect in mind the present article is understandable - we hope - for  a person who may not be familiar with operator spaces and completely bounded linear or multilinear maps.

  In this article  we are mostly interested in the matrix $X$ and not so much in the operators $F_X, B_X, S_X, T_X$, so  we will introduce 6 norms on $M_{(m,n)}(\bc),$ as the norms or completely bounded norms of the operators mentioned.  We will then quote one of the main results from \cite{C2}, where the norms are described as optimal solutions to some factorization problems. 
 
 \begin{definition}
 Let $X $ be a complex $m \times n $ matrix then 
 \begin{itemize} 
 \item[(i)] $\|X\|_F := \|F_X\|,$ named the F-norm. 
 \item[(ii)] $\|X\|_{cbF} := \|F_X\|_{cb},$ named the cbF-norm, or the completely bounded F-norm.
 \item[(iii)] $\|X\|_B := \|B_X\|,$ named the B-norm or the bilinear form norm.
 \item[(iv)] $\|X\|_{cbB}  := \|B_X\|_{cb},$
named the cbB-norm or the completely bounded  bilinear form norm. 
\item[(v)] $\|X\|_S := \|S_X\| = \|S_X\|_{cb},  $ named the S-norm or the  Schur  norm.
 \item[(vi)] $\|X\|_T  := \|T_X\|= \|T_X\|_{cb},$ named the T-norm, or the bilinear Schur norm.
 \end{itemize}
\end{definition}
 We remind you that for a vector $\xi$ in $\bc^n$ the expression $\Delta(\xi)$ symbolizes the diagonal $n \times n $ matrix with the  entries $\xi_j$ placed in the canonical way. Further we recall that for a complex matrix  $X$ the expression $\|X\|_\infty$ denotes its operator norm and the expression $\|X\|_2$ its Hilbert-Schmidt norm. The column norm $\|X\|_c$  of $X$ is defined as the maximal norm attained by the columns in $X.$

 The results on factorizations from Section 2 of \cite{C3} are summarized in the statement 1.7 from the introduction of \cite{C3}. The theorem just below reproduces the statement 1.7, except for some slightly changed formulations of the items (iii) and (iv) below. The proof of the changed versions of items (iii) and (iv)  follow right  after the formulation of the theorem. In the presentation below  , the aspect, that the norms are minimal solutions to some factorization problems, is emphasized by some displayed formulae. 
 
 \begin{theorem} \label{CbFac}
 Let $X $ be a complex $m \times n $ matrix then. 
 \begin{itemize} 
 \item[(i)] There exists a matrix $A$ in $M_{(m,n)}( \bc)$ and a vector $\xi$ in $\bc^n$ such that all $\xi_j \geq 0,$ $\|\xi\|_2 = 1 $, $\|A\|_\infty = \|X\|_{cbF}, $ the support of $A$ is dominated by that of $\Delta(\xi),$ $X = A\Delta(\xi)$ and \begin{align*} \|X\|_{cbF} = \min\{\|\hat A \|_\infty \|\Delta(\hat \xi)\|_2 \, : \, &\hat A \in M_{(m,n)}(\bc), \hat \xi \in \bc^n, \\& \text{ such that }    X = \hat A \Delta(\hat \xi)\,\}.\end{align*}  
 \item[(ii)]  There exists a matrix $B$ in $M_{(m,n)}( \bc),$ a vector $\eta$ in $\bc^m$ and a vector $\xi$ in $\bc^n$ such that all $\eta_i \geq 0,$ $\|\eta\|_2 =1, $ all $\xi_j \geq 0,$ $\|\xi\|_2 = 1 $, $\|B\|_\infty = \|X\|_{cbB} ,$  the support of $B$ is dominated by that of $\Delta(\xi),$ the range of $B$ is contained in the support of $\Delta(\eta),$  $X =\Delta(\eta) B\Delta(\xi)$  and \begin{align*}  \|X\|_{cbB} = \min\{ \|\hat \eta\|_2 \|\hat B\|_\infty\|\hat \xi\|_2 \, : \,&  \hat \eta \in \bc^m , \hat B \in M_{(m,n)}(\bc), \hat \xi \in \bc^n, \\ & \text{ such that } X = \Delta(\hat \eta) \hat B \Delta( \hat\xi)\}.
 \end{align*} 
 \item[(iii)]Let $r$ denote the rank of $X$  and $k $ a natural number such that $k \geq r.$  There exists a complex $k \times m$ matrix $L$ and a complex  $k \times n $ matrix $R$ such that $ X=  L^*R,$ $\|X\|_S = \|L\|_c\|R\|_c,$ the range projections for $L$ and $R$ are the same projection and
  \begin{align*} \|X\|_S = \min\{\|\hat L\|_c \|\hat R\|_c \, :&\, \exists k \in \bn,  \hat L,  \in M_{(k,m)}(\bc),  \hat R \in M_{(k,n)}(\bc)\\ & \text{ such that }  X = \hat L^* \hat R\}.\end{align*} 
 \item[(iv)] Let $r$ denote the rank of $X$ and $k$ a natural number such that $k \geq r.$
 There exists a vector $\eta$ in $ \bc^m,$ a complex $k \times m$ matrix $L$ and a complex  $k \times n$ matrix $R$ such that all $\eta_i \geq 0, $ 
 $\|\eta\|_2 =1,$ $X = \Delta(\eta )L^*R, $  $\|X\|_T  = \|L\|_c\|R\|_c , $  the range projections of $R,$ $L$ equal each other and \begin{align*}
 \|X\|_T = \min\{ \|\hat \eta\|_2 \|\hat L\|_c \|\hat R\|_c \,& : \, \hat \eta \in \bc^m, \exists k \in \bn, \hat L \in M_{(k,m)}(\bc) ,\\ &  \hat R \in M_{(k,n)}(\bc) \text{ such that } X = \Delta(\hat \eta) \hat L^* \hat R
\}. \end{align*}  
 \end{itemize}
\end{theorem}

\begin{proof}
The only real differences between this theorem and the statement 1.7 of \cite{C3} can be found  in the formulations of the items (iii) and (iv). In \cite{C3} Theorem 2.7 the  item (iii) is proven for $k = r,$ and since $\bc^r$ naturally embeds isometrically into $\bc^k$ for any $k \geq r$ the new formulation of item (iii) follows. The changed formulation of item(iv) may be verified by the same arguments used in the case of item (iii).
\end{proof}  The point in the change of the formulation of the items (iii) and (iv) comes from the fact, which  we want to underline, that it is always  possible to change a given factorization $X = L^*R$  with $L$ in  $M_{(k,m)}(\bc)$ and $R$ in  $M_{(k,n)}(\bc)$ into one in the form  $X = L^*_1R_1$ with $L_1 $ in  $M_{(k,m)}(\bc)$ and $R_1$ in  $M_{(k,n)}(\bc)$ such that both matrices have the same range projection. This formulation turns out to be easier to apply than the formulation where we  demand that the factors $L$ and $R$ both have $\bc^r $ as their  common range.   We formulate this as a lemma. 
  
\begin{lemma} Let $ X$ be in $M_{(m,n)}(\bc), $ $L$ in  $M_{(k,m)}(\bc)$ and $R$ in  \newline $M_{(k,n)}(\bc)$ such that $X = L^*R.$  There exists $L_1$ in  $M_{(k,m)}(\bc)$  with $\|L_1\|_c \leq \|L\|_c ,$  $\|L_1\|_2\leq \|L\|_2,$  $\|L_1\|_\infty \leq |L\|_\infty  ,$  and $R_1$ in  $M_{(k,n)}(\bc)$  with $\|R_1\|_c \leq \|R\|_c,$ $\|R_1\|_2\leq \|R\|_2,$ $\|R_1\|_\infty  \leq \|R\|_\infty $such that $X = L^*_1R_1$ and  both matrices have the same range projection.
\end{lemma} 
  \begin{proof} 
  Let $Q$ denote the range projection of $L$ and replace $ R$ by $R_1$ defined by $R_1 := QR.$ The define $Q_1$ as the range projection of $R_1$ and $L_1 := Q_1L.$ Since $Q_1 \leq Q$ we have that the range projection of $L_1 $ equals $Q_1.$  Then the range projections of both $L_1$ and $R_1$ equals $Q_1,$ and the lemma follows by well known results on the behaviour of the norms. \end{proof}

The main result, Theorem \ref{unique},  in this article states that the optimal factorizations in the items (i),  and (ii) above are uniquely determined by the optimality. The  factorizations given in (iii) and (iv) are never uniquely determined, but we define two related but more complicated factorizations in Definition \ref{AdvancedDef}. We show that the new  factorizations in the  bilinear Schur case, corresponding to item (iv), is uniquely determined. We show by a simple example that  our advanced factorization for Schur multipliers is not always uniquely determined, but we have found an extra condition under which the uniqueness is proven to hold. It turns out that the advanced Schur factorization is uniquely determined for unitaries,  for  positive matrices and there is a hope, that this type of Schur factorization may be uniquely determined for self-adjoint matrices ? This makes it reasonable for us to introduce names for  the factorizations we already have met in  the items (i)..(iv).  

\begin{definition} \label{definitions} Let $X$ be a complex $m \times n $ matrix 
\begin{itemize}
\item[(i)] A factorization $X = A\Delta(\xi) $ with the properties given in item (i) from Theorem \ref{CbFac} is called a {\em cb operator factorization} of $X.$ 
\item[(ii)] A factorization $X = \Delta(\eta) B\Delta(\xi) $ with the properties given in item (ii) from Theorem \ref{CbFac} is called a {\em cb bilinear form factorization} of $X.$ 
\item[(iii)] A factorization $X = L^*R $ with the properties given in item (iii) from Theorem \ref{CbFac} is called an {\em elementary Schur factorization} of $X.$ 
\item[(iv)] A factorization $X =\Delta(\eta)  L^*R $ with the properties given in item (iv) from Theorem \ref{CbFac} is called an {\em elementary  bilinear Schur  factorization} of $X.$ 

\end{itemize}
\end{definition} 

If $X$ is self-adjoint it seems natural to ask for a self-adjoint, or more precisely a symmetric,  factorization in the items (ii) and (iii) above, and it turns out that this is possible in both cases. 

With respect to item (ii) this follows from Theorem \ref{unique}  item (ii), where it is shown that the cb bilinear form factorization is unique. This is stated as Corollary \ref{SA}. 

In the case of the elementary Schur factorization of a self-adjoint matrix $X$ we show in the following proposition, that $X$ has a factorization $X = \|X\|_S GSG$ such that $G$ is positive with $\|G\|_c \leq 1$ and $S $ is self-adjoint partial isometry with range space equal to the range of $G.$  Then recall that $\|SG\|_c \leq \|G\|_c$ and it follows that $X = (\|X\|_S^{(1/2)}G)^*(\|X\|_S^{(1/2)}SG)$ is an elementary  elementary Schur factorization.

\begin{proposition} \label{SaElmSchur}
Let $X$ be a self-adjoint matrix in $M_n(\bc)$ then there exists a positive matrix $G$ and a self-adjoint partial isometry  $S$ in $M_n(\bc)$ such that $\|G\|_c \leq 1, $ $S^2 $ equals the support projection of $G$ and $X =\|X\|_S GSG.$

If $X$ is positive with range projection $P$  then the  matrix $G$ defined by  $G:= (1/\|X\|_S)^{(1/2)}X^{(1/2)} $ is positive with $\|G\|_c = 1,$ has range projection equal to $P$ and   $X = \|X\|_S GPG.$  \end{proposition} 

\begin{proof}
Let us assume that $\|X\|_S =1,$ and let an elementary  Schur factorization be given as $X = L^*R$ with matrices in $M_n(\bc)$ and  $\|L\|_c =  \|R\|_c =1,$ then since $X$ is self-adjoint  we can write $X$ as a product of block matrices, with matrices in $M_n(\bc).$ 
\begin{equation}
X = \begin{pmatrix}
2^{(-1/2)} L^* & 2^{(-1/2)}R^* 
\end{pmatrix} \begin{pmatrix}
0 & I_n\\ I_n & 0
\end{pmatrix} \begin{pmatrix}
2^{(-1/2)} L \\  2^{(-1/2)}R 
\end{pmatrix}.
\end{equation}
The polar decomposition of the  right hand column  factor gives a positive $n \times n$ matrix $F$ such that $F^2 = (1/2) L^*L + (1/2)R^*R$  and a pair $A, B$ of $n \times n$ matrices such that $AF = 2^{-(1/2)} L ,$ $BF = 2^{-(1/2)} R $ and $A^*A + B^*B $ is the range  projection of $F.$ It should be remarked that  diag$(F^2) \leq I_n,$ since $\|L\|_c = \|R\|_c =1, $ and then $\|F\|_c \leq 1.$   We can now define a self-adjoint contraction $T$ by the equation

\begin{equation}
T: = \begin{pmatrix}
 A^* & B^* 
\end{pmatrix} \begin{pmatrix}
0 & I_n\\ I_n & 0
\end{pmatrix} \begin{pmatrix}
A \\  B 
\end{pmatrix},
\end{equation}
and then  $X = FTF,$ so for $\hat L := F $ and $\hat R:= TF$ we have $X=  \hat L^* \hat R$ with $\|L\|_c \leq 1$ and $\|R\|_c \leq 1.$ To obtain the desired self-adjoint decomposition we let $P$ denote the range and support projection of $T$ and note, that by construction $P$ is dominated by the support projection of $F.$ Let $T = S_0|T|$ denote the polar decomposition of $T,$ then $S_0 $ is a self-adjoint  partial isometry which commutes with $|T|$ and satisfies $S^2_0 =P. $ Recall that the range of $P$ is contained in the range of $F$ so the range of $|T|^{(1/2)}F $ is the range of $P.$ Hence for the polar decomposition of $VG$ of $|T|^{(1/2)}F$
 we get that $V$ is a partial isometry from the range of $G$ onto the range of $P,$ so  the matrix $S$ defined by $S:= V^*S_0V $ is a self-adjoint partial isometry from the range of $G$ onto the range of $G.$ The matrix  $ G$ is a positive matrix such that $G^2 \leq F^2$ and then $\|G\|_c \leq 1. $  By construction we have $X = GSG$ and the first part of the proposition follows.             

The second part follows from Schur's result \cite{Sc}, which tells that for $X$ positive $S_X$ is a positive mapping and then $\|S\|_X = \|\mathrm{diag}(X)\| = \|X^{(1/2)}\|_c^2,$ and the proposition follows. 
\end{proof}
 
 It is well known that Grothendieck's  inequalities imply factorization results for matrices such that the so called {\em little inequality} in   \cite{Pi3}  gives   factorization of the same nature as that of item (i) above, see \cite{Pi3} Th. 5.2.     {\em Grothendiek's inequality } in \cite{Pi3}  gives a factorization of the same nature as that in item (ii), see \cite{Pi3} Th. 2.1. This is reproduced  in the first part of \cite{C3} in the theorems 1.1 and 1.2 . Pisier's survey article \cite{Pi3} contains definitions and properties of the  Grothendieck's constants $k_G^\bc$ and $K_G^\bc$ and much more. Based on the results of \cite{C3} we have chosen to formulate our interpretation of Grothendieck's inequalities in the following way.

 \begin{theorem} \label{Gr1}
 Let $X$ be a complex  $m \times n $ matrix, then
 \begin{itemize}
 \item[(i)]  $\|X\|_{cbF} \leq \sqrt{k_G^\bc} \|X\|_F,$ and $\sqrt{k_G^\bc}$ is the best possible constant for this inequality,
 \item[(ii)] $\|X\|_{cbB} \leq K_G^\bc\|X\|_B,$ and $K_G^\bc$ is the best possible constant for this inequality. 
 \end{itemize}
 
\end{theorem}   

There is a natural inner product on $M_{(m,n)}(\bc),$ which considers \newline $M_{(m,n)}(\bc)$ to be a copy of $\bc^{mn},$  and it may be expressed via  the standard trace Tr$_n$ on  $M_n(\bc),$ as follows   
\begin{equation} \label{InProd}
\forall X, Y \in M_{(m,n)}(\bc): \quad \langle X, Y \rangle := \mathrm{Tr}_n(Y^*X).
\end{equation}

 The other result from \cite{C3}, which we will use heavily in the arguments to come,  states that this inner product  gives us a way to show that the conjugate dual to the normed space $(M_{(m,n)}(\bc), \|.\|_{cbF})$ is isometrically isomorphic to   $(M_{(m,n)}(\bc), \|.\|_{T}),$ and similarly the conjugate dual of  $(M_{(m,n)}(\bc), \|.\|_{cbB})$ equals  $(M_{(m,n)}(\bc), \|.\|_S).$  We want to formulate this using the polar operation $\cas \to \cas^\circ$ defined on subsets of  $M_{(m,n)}(\bc)$ via the inner product as 
 \begin{equation} \label{polar}
 \forall \cas    \subseteq  M_{(m,n)}(\bc): \quad \cas^\circ:= \{Y \in  M_{(m,n)}(\bc)\, : \, \forall X \in \cas:\, \, \, |\langle X, Y\rangle | \leq 1 \, \} 
\end{equation}  Let Ball$_{cbF}$, Ball$_{cbB}$, Ball$_S$ and Ball$_{T}$ denote the unit balls in $M_{(m,n)}(\bc)$ equipped with the norms with the same suffix, then the equations (3.5) and (3.6) of \cite{C3} give the following theorem.

\begin{theorem} \label{Polar} All the balls $\mathrm{Ball}_{cbF}, \mathrm{Ball}_{cbB}, \mathrm{Ball}_S$ and $ \mathrm{Ball}_{T} $ equal their bi-polars and 
\begin{itemize}
\item[(i)] $(\mathrm{Ball}_{cbF})^\circ = \mathrm{Ball}_{T}.$
\item[(i)]$ (\mathrm{Ball}_{cbB})^\circ = \mathrm{Ball}_S.$
\end{itemize}
\end{theorem}

\section{Uniqueness of some factorizations} 

We start right away by stating and proving the uniqueness results we have for the {\em cb operator factorization } and the {\em cb bilinear form factorization} given in the items (i) and (ii) of  Definition \ref{definitions}. Then we will introduce some slightly more complicated Schur factorizations corresponding to the items (iii) and (iv) and show that the new version of item (iv) is unique and the advanced Schur factorization in item (iii) is unique  for unitaries, for positive matrices and in general under an extra condition.   

\begin{theorem} \label{unique} 
Let $ X$ be a complex $m \times n$ matrix 
\begin{itemize}
\item[(i)] A cb operator factorization of $X$ is uniquely determined.
\item[(ii)] A cb bilinear form factorization of $X$ is uniquely determined. 
\end{itemize} 
\end{theorem} 

\begin{proof}
We give the proofs according to their numbering, so assume that $\|X\|_{cbF} = 1 $ and that a completely bounded operator factorization of $X$ is given by $X = A\Delta(\xi).$ First we remark, that we may assume that no column in $X$ vanishes, since the vanishing of a column will imply that the corresponding vector coordinate $\xi_j $ will vanish, because in the opposite it will be possible to construct a   factorization  with a product of the norms lower than  the cbF-norm, and this contradicts item (i) of Theorem  \ref{CbFac}.
 Then a deletion of a vanishing column will not change $\xi,$ and we may as well delete that column.  With the assumption that no column in $X$ vanishes, the given factorization of $X$ implies that no $\xi_j$ vanishes. 
 By Theorem \ref{Polar} item (i) there exists a complex $m \times n $ matrix $Y$ with $\|Y\|_T =1 $ such that Tr$_n( Y^*X) = 1, $ and $Y $ has an elementary  
bilinear Schur factorization $Y = \Delta(\eta) L^*R.$ It is elementary to compute a couple of Hilbert-Schmidt norms and  see that  $\|R\Delta(\xi)\|_2 \leq  \|R\|_c\|\xi\|_2 \leq 1$ and similarly $\|L\Delta(\eta )\|_2 \leq 1 . $ Hence 
$\|L\Delta(\eta)A\|_2 \leq \|A\|_\infty \leq 1, $ and we may compute \begin{align} \label{C=S1} 
1 & = \mathrm{Tr}_n (Y^*X)\\ \notag &= \mathrm{Tr}_n( R^*L\Delta(\eta)A\Delta(\xi))\\ \notag
&= \mathrm{Tr}_n\big((\Delta(\xi)R^*)(L\Delta(\eta)A)\big)\\ \notag &= \langle L\Delta(\eta)A, R\Delta(\xi) \rangle  \\ \notag &\leq \|R\Delta(\xi)\|_2 \|L\Delta(\eta)A\|_2 \\ \notag &\leq 1. 
\end{align}
This implies that $\|L\Delta(\eta)\|_2 = \|R\Delta(\xi)\|_2 = 1$ and $L\Delta(\eta) A = R \Delta(\xi).$ Since we already know that  $\|R\|_c  = \|\xi\|_2 = 1,$ and all $\xi_j > 0,$ the equality  $\|R\Delta(\xi)\|_2 = 1$ and  an elementary calculation of $\|R\Delta(\xi)\|_2^2$  show  that all columns in $R$ must be of unit length and in particular diag$(R^*R) = I_n.$ We may then combine these observations as follows  \begin{equation}
Y^*X =R^*(L\Delta(\eta)A)\Delta(\xi) = R^*R\Delta(\xi)^2.
\end{equation}
Since diag$(R^*R) = I_n$ we find that 
\begin{equation}
\Delta(\xi)^2 =\mathrm{diag}(Y^*X),
\end{equation}
Since each $\xi_j > 0, $ we get that
$\xi $ and then also $A$ is uniquely determined and the proof of item (i) follows.

The proof of item (ii) is quite similar to the one just presented, so we suppose that $X$ is a complex   $m \times n $ matrix with a completely bounded bilinear form factorization $X = \Delta(\eta)B \Delta( \xi), $ and we also do  assume that no column and no row in $X$ vanishes. As above we assume further that $\|X\|_{cbB} = 1$ and that $Y $ is a complex $m \times n$ matrix with $\|Y\|_S = 1$
and an  elementary Schur factorization $Y = L^*R$ with $\|R\|_c = \|L\|_c =1,$ such that Tr$_n(Y^*X) = 1.$ Then we may follow the computations in (\ref{C=S1}) to see that
$$\langle R\Delta(\xi), L\Delta(\eta)B\rangle = \langle R\Delta(\xi)B^* , L\Delta(\eta)\rangle =1.$$ 
Since $\|B\|_\infty \leq 1 ,$ we get 
 $\|L\Delta(\eta)\|_2 = 1, $ $\|R\Delta(\xi)\|_2 =1, $ \newline             $\|L\Delta(\eta)B\|_2 = 1,  \, \, \|R\Delta(\xi)B^*\|_ 2 = 1,$ so  \begin{equation} \label{CS2}
\text{(i)}\quad  L\Delta(\eta) B = R\Delta(\xi) \text{ and (ii)} \quad L\Delta(\eta) = R \Delta(\xi) B^*.
\end{equation}
Since no row and no column in $X$ vanishes, we see that  all the $\xi_j >0$ and all the $\eta_i > 0.$ This implies as above that all columns in both $L$ and $R$ have unit length so 
\begin{equation} \label{diagonals} 
\mathrm{diag}(L^*L) = I_m \quad \text{ and } \quad \mathrm{diag}(R^*R ) = I_n.
\end{equation}
Based on the equation ( \ref{CS2}) we get
that for any $Y$ with $\|Y\|_S = 1 $ and Tr$_n(Y^*X) = 1$ we have 
\begin{align}
Y^*X &= R^*[L \Delta(\eta)B]\Delta(\xi) = (R^*R)\Delta(\xi)^2.\\
YX^* &= L^*[R\Delta(\xi)B^*]\Delta(\eta) = (L^*L) \Delta(\eta)^2.
\end{align}
By (\ref{diagonals}), and the inequalities $\xi_j >0,$   $\eta_i >0,$ we find that in the first place  $\xi$ and $\eta$ are uniquely determined and then also $B$ is uniquely determined via the equations  \begin{equation} \label{vectors} 
\Delta(\xi)^2 = \mathrm{diag}(Y^*X) \text{ and } \Delta( \eta)^2 = \mathrm{diag}(YX^*).
\end{equation}
The theorem follows. 
\end{proof}
There is a corollary to the uniqueness result of  Theorem \ref{unique} item (ii),  which is quite easy to obtain but still valuable. 

\begin{corollary} \label{SA} 
Let $X $ be a self-adjoint matrix in $M_n(\bc)$ with cb-bilinear form factorization $\Delta(\eta) B \Delta(\xi) $ then $\xi = \eta $ and $B = B^*.$ If $X$ is positive, so is $B.$ 
\end{corollary} 
\medskip
 Suppose a complex  matrix $X$ is given as a product $X=L^*R$ with $L $ in $M_{(k,m)}(\bc)$ and $R$  in $M_{(k,n)}(\bc).$ Then for any unitary $U$ in $M_k(\bc)$ we get a  factorization $X =(UL)^*(UR)$ with $\|UL\|_c = \|L\|_c $ and $\|UR\|_c= \|R\|_c.$ If the ranges of $L$ and $R$ agree, so do the ranges of $UL$ and $UR,$ of course. This indicates that a more elaborated point of view on the Schur factorization has to be used in order to get some sort of uniqueness, and it also suggests that we should look for yet another application of the polar decomposition. If an elementary Schur factorization of $X$ with $\|X\|_S = 1$   is given as $X= L^*R$ with $\|L\|_c = \|R\|_c = 1, $ then we should look for the polar decompositions of $L$ and $R,$ say $L = V_L|L|$ and $R = V_R|R| $ with $V_L$ a partial isometry from the support of $L$ onto the range of $L$ and similarly for $V_R.$ 
 Since the ranges of $L$ and $R$ are the same space, we can define a partial isometry $ W: = V_R^*V_L $ in $M_{(m,n)}(\bc)$ with support equal to the support of $|L|$ and range equal to the support of $|R| $  and we have $X = |R| W |L|,$ such that $|L| \geq 0 ,  $ diag$(|L|^2) \leq I_n,$ $|R| \geq 0 ,$ diag$(|R|^2 ) \leq I_m$ and $W$ is a partial isometry from the support of $|L|$ onto the support of $|R| .$  

 A good part  of this was discussed previously in the article \cite{C2}, where we also discussed block matrices. Here we summarize this in the item (i) in the definition below, but we will first look at the elementary bilinear Schur factorization, which is clearly not unique for the same reasons as given above. 
 First we will show that without loss of generality we may,  in a discussion of the bilinear Schur factorization, as well assume that  all rows in $X$ are non trivial. If row number $i$ in $X$ is trivial and  $X = \Delta(\eta) L^*R$ is an optimal elementary bilinear Schur factorization, then $\eta_i = 0,$ because a valid inequality such as $\eta_i >0$ will spoil the optimality condition in the factorization.  We may then just as well take the $i$'th row of $X$ out before we start factorizing. We then observe that in  an elementary bilinear Schur factorization of an $X,$  such as $X = \Delta(\eta )L^*R,$ with $\|\eta\|_2 = \|X\|_T,$ and $\|L\|_c = \|R\|_c = 1, $ where all rows in $X$  are non trivial, we have all $\eta_i >0$ and
 $\|L^*R\|_S = 1 .$ Based on the arguments just above, we may replace $L^*R$ by a product $FW_0G$ such that $F$ is positive in $M_m(\bc)$ with diag$(F^2) \leq I_m,$ $G$ is positive in $M_n(\bc) $ with diag$(G^2) \leq I_n$ and $W_0$ is a partial isometry in $M_{(m,n)}(\bc)$ with support equal to the range of $G$ and range equal to the support of $F.$ 
 Then $X = \Delta(\eta)FW_0G.$ If just one row, say row number 1,   in $F$ has norm less than $1,$  then there exists a positive real $t$ with $0 < t < 1 $ such that if the first row of $F$ is multiplied by $t^{-1},$ then the row norm of the new matrix, say $\hat F$  is still at most 1, and then  the Hilbert space norm of the vector $\hat\eta,$ obtained from $\eta$ by replacing the first entry $\eta_1$ by $t\eta_1,$ is smaller than $\|X\|_T$  since all $\eta_i >0.$ 
 Then $X = \Delta(\hat \eta) \hat F ( W_0 G)$ and by Theorem \ref{CbFac} item (iv)  $\|X\|_T \leq \|\hat\eta\|_2 < \|X\|_T.$ Hence all rows in $F$ have norm 1 and the Hilbert-Schmidt norm $\|\Delta(\eta)F\|_2$ equals $\|X\|_T.$  
  Let $F\Delta(\eta) $ have the polar decomposition $VT$ such that $T$ is a positive $m \times m $ matrix and $V$ is a partial isometry from the support of $T$ onto the range of $F\Delta(\eta)$ which in turn equals the range and also support space  of $F$ since $\Delta(\eta) $ is invertible and $F$ is positive. 
 Let $W := V^*W_0,$ then  $W$ is a partial isometry from the range space of $G$  onto the support space of the positive matrix $T, $ and we have obtained a factorization of $X$ as $X = TWG.$ 
  Here $T$ is positive in $M_m(\bc)$ with $\|T\|_2=\|X\|_T, $ $G$ is positive in $M_n(\bc)$ with $\|G\|_c = 1 $ and $W$ is an isometry in $M_{(m,n)}(\bc) $ with support equal to the range of $G$ and range equal to the support of $T.$ 
 
\begin{definition} \label{AdvancedDef}
Let $X$ be a complex $m \times n $ matrix. 
\begin{itemize}
\item[(i)] A factorization $X = \|X\|_SFWG, $ such that $F$ is a positive $m \times m $ matrix with diag$(F^2) \leq I_m,$  $G$ is a positive  $n \times n $ matrix with diag$(G^2) \leq I_n$ and $W$ is an $m \times n $ partial isometry   from the range of $G$ onto the range of $F$ is called a Schur factorization of $X.$ 
\item[(ii)] A factorization $X = TWG, $ such that $T$ is a positive   $m \times m $ matrix with $\|T\|_2 = \|X\|_T ,$   $G$ is a positive  $n \times n $ matrix with diag$(G^2) \leq I_n$ and $W$ is a partial isometry from the range space of $G$ onto the support space of $T$   is called a  bilinear Schur factorization of $X.$ 
\end{itemize}
\end{definition}

The lines before the definition above yield the first part of the following proposition. 

\begin{proposition} 
Let $X$ be a complex $m \times n $ matrix, then $X$ has both a Schur factorization and a bilinear Schur factorization. If $X$ is self-adjoint it is possible to obtain a Schur factorization of the form $X=\|X\|_S GSG$ with $S$ a self-adjoint partial isometry. 
\begin{itemize}
\item[(i)] If $X= FCG$ with   $ F $ positive in $M_m(\bc),$   $C$ in $M_{(m,n)}(\bc)$  and $G$ positive in $M_n(\bc) $ then $$\|X\|_S \leq\|\mathrm{diag}(F^2)\|^{(1/2)} \|C\|_\infty\|\mathrm{diag}(G^2)\|^{(1/2)}.$$ 
\item[(ii)] If $X = TCG$ with $T$ positive  in $M_m(\bc),$  $C$ in $M_{(m,n)}(\bc) $ and $G$ positive in $M_n(\bc)$ then $$\|X\|_T \leq \|T\|_2\|\|C\|_\infty \|\mathrm{diag}(G^2)\|^{(1/2)}.$$  
\end{itemize} 

\end{proposition} 

\begin{proof}
The existence of the factorizations were demonstrated in front of  the definition above. The special result for Schur factorizations of self-adjoint matrices follows from Proposition \ref{SaElmSchur}.

With respect to item (i) you may write $X = F^*(CG)$ and by item (iii) of Theorem \ref{CbFac} we have $$\|X\|_S \leq \|F\|_c\|CG\|_c \leq \|\mathrm{diag}(F^2)\|^{(1/2)} \|C\|_\infty\|\mathrm{diag}(G^2)\|^{(1/2)},$$ and item (i) follows.

For item (ii) we define  a positive  vector $\eta$ 
in $\bc^m$ by \newline
$\eta_i := (\sum_j |T_{(i,j)}|^2)^{(1/2) }  $ and a matrix $L$ in $M_m(\bc)$ by the equation $T = \Delta(\eta)L^*,$ and the demand that if row number $k$ in $T$ vanishes, so does column  number $k$ of $L.$  Then $\|\eta\|_2 = \|T\|_2, $ $L_{(s,t)} = \overline{ T_{(t,s)}}/ \eta_t$ and for any $t$ we have $\sum_s |L_{(s,t)}|^2 = 1,$
 so  all columns in $L$ have norm $1$  and  $\|L\|_c =1.$ If we define $R := CG$ then as in the proof of item (i), we have $\|R\|_c \leq \|C\|_\infty\| \mathrm{diag}(G^2)\|^{(1/2)}. $ Now $X = \Delta(\eta)L^*R$ with  $\|\eta\|_2 \|L\|_c \|R\|_c \leq \|T\|_2 \|C_\infty\|\| \mathrm{diag}(G^2)\|^{(1/2)}, $ and the proposition follows from item (iv) of Theorem \ref{CbFac}.

\end{proof}

With respect to the uniqueness of these factorizations we show below that the bilinear Schur factorization of a matrix is uniquely determined. The Schur factorization can not be unique in many cases as the following example will show.

\begin{example} In this example $i$ denotes the imaginary complex unit.  Let $X$  in $M_2(\bc)$ given by $$X = \begin{pmatrix}
1 & 0\\0 & (i/4)
\end{pmatrix},$$ then $\|X\|_S = 1,$ and  we may give several possible Schur factorizations such as 
$$X = \begin{pmatrix}
1 & 0\\0 & (1/2)
\end{pmatrix}\begin{pmatrix} 1 & 0\\0  & i
\end{pmatrix}\begin{pmatrix} 1 & 0\\0 & (1/2)
\end{pmatrix} =\begin{pmatrix}
1 & 0\\0 & 1
\end{pmatrix}\begin{pmatrix} 1 & 0\\0 & i
\end{pmatrix}\begin{pmatrix} 1 & 0\\0 & (1/4)
\end{pmatrix}.$$
The things which makes this possible is that a sub-matrix of $X$ has Schur multiplier equal to the Schur norm of $X,$ and  that this $X$ is not self-adjpint.  We will show in Propositon  \ref{SameNorm} that a matrix $X,$  which has the property that any sub matrix obtained by deleting a column has smaller Schur norm than $\|X\|_S,$ has a unique Schur factorization. With respect to the self-adjoint case, the example above can easily be modified to give two different Schur factorizations of a self-adjoint matrix, but these will not be symmetric  in the sense that the left hand positive factor equals the right hand positive  factor. We have not been able to construct two different and  symmetric Schur factorizations of a self-adjoint matrix, so may be there is only one symmetric Schur factorization of a self-adjoint matrix ? 

It is easy to see that a positive matrix $X$ has the unique Schur factorization given as $X = \|X\|_S(\|X\|_S^{-(1/2)}X^{(1/2)})P(\|X\|_S^{-(1/2)}X^{(1/2)}),$ where $P$ is the range projection of $X.$ 

We show  in Theorem \ref{unitaries}  that a unitary matrix in $M_n(\bc)$ has the unique Schhur factorization given by $U = (I_n) U (I_n).$ \end{example}  

The remarks at the end of the previous example  indicates that the  factorization we have chosen to give the name {\em Schur factorization }  may be reasonable,  but there are other possibilities for a good  definition of what a Schur factorization can  be. The following proposition could be the basis for  an alternative definition of a Schur factorization, but it seems to be less detailed and  not so closely related to the bilinear Schur factorization, which by Theorem \ref{UniBilSchur} is uniquely determined. 

\begin{proposition} Let $X$ a be a complex $m \times n  $ matrix for which $\|X\|_S = 1,$ then $X$ has a factorization $ X = FCG$ such that $\|C\|_\infty =1, $ $F$ and $G$ are positive with  $\mathrm{diag}(F^2) = I_m$ and $\mathrm{diag}(G^2) = I_n.$ 
\end{proposition}

\begin{proof}
Let $X = L^*R$ be an elementary Schur factorization of $X$ with $\|L\|_c = \|R\|_c = 1,$  and let $k$ denote the number of rows in $R$ and $L.$ We define two new matrices $\hat L$ and $\hat R$ by adding two more rows to both $L$ and $R.$ The first extra row in $\hat L$ consists of the numbers $(1 -\|L_j\|^2)^{(1/2)},$ where $L_j$ denotes the $j'$th column in $L,$  and the next row vanishes. 
Then all columns in $\hat L$ are unit vectors. We do nearly the same in $\hat R$ except that here the first extra row vanishes and the last one supplements $R$ such that all columns in $\hat R$ are unit vectors. Clearly $\hat R^* \hat L = X $ and for the polar decompositions $\hat L = V|\hat L| $ and $\hat R = W |\hat R| ,$ we get the promised  factorization as $ X = |\hat R|( W^*V) |\hat L|.$ 
    \end{proof}

\begin{proposition} \label{SameNorm}
Let $X $ be a complex $m \times n $ matrix. If $X$ has the property that  the Schur norm of any sub-matrix obtained by deleting a column  is  smaller than the Schur norm of $X$, then the Schur factorization is unique. 

If the Schur norm of any sub-matrix obtained by deleting a row  is  smaller than the Schur norm of $X$, then the Schur factorization is unique. 
\end{proposition}

\begin{proof}
Suppose $\|X\|_ S = 1 $ and let a Schur factorization $X = FWG$  be given. Suppose $Y$ is  a complex $m \times n $ matrix such that $\|Y\|_{cbB} =1 $ and Tr$_n(Y^*X) = 1.$ The cb bilinear form  factorization of $Y$ is given as $Y = \Delta(\eta)B\Delta(\xi)$ with $\|B\|_\infty = 1,$ then - as in the proof of Theorem \ref{unique} -  we get $\|G \Delta(\xi)\|_2 \leq 1$ and $\|\Delta(\eta) F\|_2 \leq 1, $  and we may compute 
\begin{align}
1 & = \mathrm{Tr}_n\big(\Delta(\xi)B^* \Delta(\eta)F W G\big) \\ \notag & = \mathrm{Tr}_n \big((G\Delta(\xi))(B^*\Delta(\eta)FW)\big) \\
\notag & = \langle B^*\Delta(\eta)FW, \Delta(\xi)G \rangle \\
 \notag & \leq \| B^*\Delta(\eta)FW\|_2  \|\Delta(\xi)G\|_2  \\
 \notag & \leq 1.
\end{align}
Then we see that \begin{equation} \label{deltaG} 
\|\Delta(\xi)G\|_2 = \|B^*\Delta(\eta)FW\|_2  = 1 \text{ and } \Delta(\xi)G = B^*\Delta(\eta)FW.
\end{equation}
Since the Schur norm of a submatrix of $X$ obtained by deleting a column of $X$ is strictly smaller than the norm $\|X\|_S,$ we can see that all $\xi_j >0.$ On the other hand this and the fact that $\|\Delta(\xi)G\|_2 = 1$ implies that all rows in $G $ have unit length so diag$(G^2) = I_n.$ With this in mind we may return to (\ref{deltaG} ) and obtain \begin{equation}
\Delta(\xi)^2 G^2 = \Delta(\xi)( B^*\Delta(\eta)FW) G= Y^*X.
\end{equation}  
Since diag$(G^2) = I_n $ we get \begin{equation} \Delta(\xi)^2  = \mathrm{diag}(Y^*X). \end{equation} 
We know that all $\xi_j >0,$ so  $G^2 = \mathrm{diag}(Y^*X)^{-1}Y^*X,$ and the positive right hand factor in any Schur factorization of $X$ must be $G.$   

Let $ X = F_1W_1 G = F_2W_2 G $ be Schur factorizations of $X.$  Since the range of $G$ equals the support of both $W_1$ and $W_2,$ we get $F_1W_1 =  F_2 W_2$ and since both of these have the form of the adjoints of polar decompositions we get that $F_1 = F_2$ and $W_1 = W_2,$ so the uniqueness of the Schur factorization of $X$ is established. 

The proof of the second statement in the proposition follows by an application of the result proven to $X^*,$  and the proposition follows.
\end{proof}

\begin{theorem} \label{UniBilSchur}
Let $X$ be a complex $m \times n $ matrix then there is only one possible bilinear Schur factorization of $X.$ 
\end{theorem}

\begin{proof} Suppose that $\|X\|_T = 1,$ and let $X=TWG$ be a bilinear Schur factorization of $ X.$ Choose a complex $m \times n $matrix $Y$ such that $\|Y\|_{cbF} = 1 $ and Tr$_n(Y^*X) = 1,$ then $Y = A \Delta( \xi) $ with $\|A\|_\infty = 1 $ and $\xi$ is a positive unit vector in $\bc^n$. Remember that $\|A\Delta(\xi) G W^*\|_2  \leq 1$ and $\|T\|_2 = 1,$ then \begin{align*}
1 & = \mathrm{Tr}_n (\Delta(\xi)A^* TWG)\\& = \mathrm{Tr}_m \big((A\Delta(\xi)GW^*)^* T\big) 
\\& = \langle T, A\Delta(\xi)GW^*\rangle \\ & \leq \|T\|_2\|A\Delta(\xi)GW^*\|_2  \\ &\leq 1 .
\end{align*} 
As before, equality in the Cauchy-Schwarz inequality implies that  $T = A\Delta(\xi)GW^*$ and then  $$ T^2 = (A\Delta(\xi)GW^*)T = YX^*, $$ so $T^2$ and then, because  $T$ is positive, $T$ is uniquely determined. Since $X = TWG$ and the range of $W$ equals the support of $T$ we have that  $WG$ is uniquely determined.  This last product has the form of the polar decomposition of itself, so it is uniquely determined, and the theorem follows. 
 \end{proof}

\section{Some applications}

The uniqueness and the optimality of  some factorizations may be used to show that a given factorization has certain properties. This principle applies to unitaries and  isometries.
The principle  will also be used to make a connection between the norms $\|.\|_{cbF} $ and $\|.\|_{cbB}.$     
We end this section by showing that the cb bilinear form factorization may be obtained as the solution of    a linear program with infinitely many constraints.

\begin{theorem} \label{unitaries} 
Let $U$ be a unitary  in $M_n(\bc)$ then $\|U\|_F = \|U\|_{cbF} = \sqrt{n}, $  $\|U\|_{cbB} = n,$  $\|U\|_S = 1$ and  $\|U\|_T = \sqrt{n}.$ The F-factorization of $U$ equals $(n^{(1/2)}U)(n^{-(1/2)}I_n),$ the cb bilinear factorization  is $U = (n^{-(1/2)}I_n) (nU) (n^{-(1/2)}I_n)$  The  Schur  factorization of $U$ is uniquely determined and equals $U =I_n (U ) I_n$ and the T-factorization is $U = (U)(I_n).$  
 \end{theorem} 
 \begin{proof}
Since the columns in $U$ is a set of pairwise orthogonal unit vectors we get $\|U\|_F \leq \sqrt{n}.$ On the other hand the vector $\Omega$ in $\bc^n,$ which has all entries equal to 1, has the property that  $\|U\Omega\|_2 = \sqrt{n},$  so $\|U\|_F = \sqrt{n}.$ 

The factorization $U = (n^{(1/2)}U)(n^{-(1/2)}I_n)$ shows by Theorem \ref{CbFac} item (i) that $\|U\|_{cbF} \leq \sqrt{n} = \|U\|_F,$ so $\|U\|_{cbF} = \sqrt{n}$ too. 

Item (iv) in Theorem \ref{CbFac}  shows that 
the factorization, given for the T-norm of $U$ above, implies  that $\|U\|_T \leq \|U\|_2 = \sqrt{n}.$ On the other hand$$n = \mathrm{Tr}_n(U^*U) \leq \|U\|_{cbF}\|U\|_T \leq  n ,$$ so $\|U\|_T = \sqrt{n},$ and we have got the right factorization in the F and in the T case. 

 The factorizations of $U$ given for the cb bilinear form case and the Schur case  show by Theorem \ref{CbFac} that $\|U\|_{cbB} \leq n$ and $\|U\|_S \leq 1.$  On the other hand $$n = \mathrm{Tr}_n(U^*U) \leq \|U\|_{cbB}\|U\|_S \leq n,$$ so we have $\|U\|_{cbB} = n$ and $\|U\|_S = 1.$ The completely bounded  bilinear form factorization is then as claimed in the proposition. Let $U = FWG$ be a Schur factorization of $U,$ then $W$ must have rank $n,$ so it is  a unitary and we have $$ n = \mathrm{Tr}_n(U^*FWG ) = \langle WG, FU\rangle \leq \|WG\|_2\|FU\|_2 \leq (\sqrt{n})^2 .$$ Hence $WG = FU $ and we know  $U = FWG$ so $U = F^2U $ and $F^2 = I_n .$ Then $F=I_n$ and similarly  $E = I_n,$ so the Schur factorization of a unitary is unique as $U = (I_n)U(I_n).$  With respect to the bilinear form norm, $\|U\|_B$ norm, we know by Grothendieck's inequality that $\|U\|_B \geq \|U\|_{cbB} (K_G^\bc)^{-1} =   n(K_G^\bc)^{-1}.$ 
 \end{proof}

Our next theorem implies that there is a link between the little Grothendieck constant $k^\bc_G$ and the big one $K_G^\bc.$ We will comment on this in a forthcoming article.

\begin{theorem} \label{X*X}
Let $X$ be a complex $m \times n $ matrix then 
\begin{itemize} 
\item[(i)] $\|X\|_F^2 = \|X^*X \|_B$
\item[(ii)] $ \|X\|_{cbF}^2  = \|X^*X\|_{cbB}.$
\end{itemize}
\end{theorem}  
  
\begin{proof}
With respect to item (i), the inequality $\|X\|_F^2 \leq \|X^*X\|_B$ is obvious and the opposite inequality follows from the Cauchy-Schwarz inequality. The item (ii)  equality my be seen using the factorization results twice. Let $X = A \Delta(\xi) $ be a cb-F factorization of $X$ then $X^*X = \Delta(\xi) (A^*A) \Delta(\xi) $ and by Theorem \ref{CbFac} items (i) and  (ii) $\|X^*X\|_{cbB} \leq \|X\|_{cbF}^2.$  By Corollary \ref{SA} let  the cb bilinear factorization of $ X^*X$ be given as $X^*X  = \Delta(\xi) B \Delta(\xi),$ then there exists a partial isometry $W$ such that $X = (WB^{(1/2)} )\Delta(\xi) $ and we get from Theorem \ref{CbFac} item (i) that  $\|X\|_{cbF} \leq \|X^*X\|_{cbB}^{(1/2)}.$  
\end{proof}

Based on the uniqueness  of the cb-F factorization and  of  the cb-B factorization, the theorem above  has an immediate corollary.

\begin{corollary} \label{squares}
Let $ A, X$ be  complex $ m \times n$ matrices  and $\xi$ a unit vector in $\bc^n$ with non negative entries such that $X = A\Delta(\xi).$ The factorization $X= A \Delta(\xi)$ is the cb operator factorization of $X$ if and only if $X^*X = \Delta(\xi) (A^*A) \Delta(\xi) $ is the cb bilinear form  factorization of $X^*X.$    
\end{corollary} 

\begin{corollary} \label{isomFnorm}  Suppose $m \geq n.$
If $X$ is an isometry then
 \begin{itemize}
\item[(i)] the cb-F   factorization is given as $X = \newline  (n^{(1/2)} X)(n^{(-1/2)} I_n),$ and $\|X\|_{cbF} = \sqrt{n} = \|X\|_F,$ 
\item[(ii)] the bilinear Schur factorization is given as $X =PX,$ where $P$ is the range  projection for $X,$ and $\|X\|_T = \sqrt{n}.$  \end{itemize}
\end{corollary} 

\begin{proof}
Since $X$ is an isometry, its columns are pairwise orthogonal  unit vectors in $\bc^m,$ and it follows that $\|X\|_F = \sqrt{n}.$ For the proof of the rest of item (i) recall that $X^*X = I_n$ and then  combine Corollary \ref{squares} with Theorem  \ref{unitaries}.

For item (ii) we see that since $X$ is an isometry, the range projection $P$ must satisfy $\|P\|_2 = \sqrt{n},$ and for $X$ we have $\|X\|_c =1 $ so for $X = PX,$ we have $\|X\|_T \leq \sqrt{n}.$ On the other hand we know that $\|X\|_{cbF}  = \sqrt{n}$ and Tr$_n(X^*X) = n $ so $\|X\|_T \geq \sqrt{n} ,$ and the corollary follows.  
\end{proof}

The results on positive bilinear forms,  show that we  may compute the cb bilinear form  norm of a positive matrix, and the cb bilinear factorization as the solution of a minimization program.  

\begin{theorem}
Let $P$ be a  non zero positive complex $n \times n $ matrix then 
$$\|P\|_{cbB} = \min\{ \mathrm{Tr}_n(\Delta(\g))\, : \, \g \in \bc^n ,\g_j \geq 0, \, P \leq \Delta(\g)\, \},$$  and the $\Delta(\xi)$ in the cb bilinear  factorization of $P$ is given by an optimal solution $\g$ as  $$\Delta(\xi) = \big(\mathrm{Tr}_n(\Delta(\g))^{-(1/2)}\big)\Delta(\g)^{(1/2)}.$$
The program has a unique solution. 
\end{theorem}

\begin{proof}
Suppose $\g$ is a non negative vector in $\bc^n$ such that $P \leq \Delta(\g).$ Then Tr$_n(\Delta(\g)) > 0$ and we can define the non negative unit vector $\xi$ as in the formulation of the theorem. 
The inequality $P \leq \Delta(\g)$ implies that there exists a positive  matrix $B$ with support dominated by the support of $\Delta(\xi)$  such that $$0 \leq B \leq \mathrm{Tr}_n(\Delta(\g)) I_n\, \text{ and }  \Delta(\xi) B \Delta(\xi) = P \leq  \mathrm{Tr}_n(\Delta( \g)) \Delta(\xi)^2.$$  Then by Theorem \ref{CbFac} item (ii) we see that $\|P\|_{cbB} \leq \mathrm{Tr}_n(\Delta(\g)).$ On the other hand Corollary \ref{SA} implies that  we have the cbB factorization $P = \Delta(\xi) B \Delta(\xi ) $ with $0  \leq B$ and $\|B\|_\infty =  \|P\|_{cbB},$ so
$$P = \Delta(\xi) B \Delta(\xi ) \leq \|B\| \Delta(\xi)^2 = \|P\|_{cbB}\Delta(\xi)^2, $$  and we see that the infimum is attained, so the theorem follows.  
\end{proof}

\begin{remark}
It seems worth to remark, that the minimization task we presented above is in fact a linear program except that the number of linear constraints  is infinite. To see this we write it as a a linear program 
\begin{align*}   \text{Min: }\quad \quad \, \, \quad   \l_1 + \dots + \l_n &\\ \text{ subject to } \, \,  \forall \, j \in \{1, \dots, n\}: \l_j &\geq 0\\ 
\forall \xi \in \bc^n : |\xi_1|^2 \l_1 + \dots + |\xi_n|^2 \l_n  &\geq  \langle P \xi, \xi\rangle.
\end{align*} 
It is of course sufficient to look at unit vectors $\xi$  only, the demand that $\l_j \geq 0$ can be omitted    and much more can probably be done to obtain a better form.  In Section 6 of \cite{HJS} by  Holbrook, Johnston and Schoch they  define a linear program to evaluate the Schur norm of an $n \times n $ complex matrix, and it seems that especially the dual version of their program may be related to the methods we have been using in this article ? 
\end{remark}

\end{document}